\documentclass[11pt]{amsart}

\usepackage[marginratio=1:1]{geometry}
\usepackage{amsmath}
\newtheorem{theorem}{Theorem}[section]
\newtheorem{lemma}[theorem]{Lemma}

\newtheorem{proposition}[theorem]{Proposition}

\theoremstyle{definition}
\newtheorem{definition}[theorem]{Definition}
\newtheorem{remark}[theorem]{Remark}
\newtheorem{example}[theorem]{Example}

\newcommand{\R}{\mathbb{R}}

\newcommand{\spn}{\operatorname{span}}
\newcommand{\tr}{\operatorname{tr}}
\newcommand{\rk}{\operatorname{rk}}
\newcommand{\ad}{\operatorname{ad}}
\newcommand{\Sec}{\operatorname{Sec}}

\newcommand{\VV}{\mathcal{V}}
\newcommand{\DD}{\mathcal{D}}
\newcommand{\HH}{\mathcal{H}}

\newcommand{\la}{\langle}
\newcommand{\ra}{\rangle}
\newcommand{\QED}{\hfill $\Box$}

\newcommand{\e}{\varepsilon}

\title{Conjugate points of dynamic pairs and control systems}
\author{B. Jakubczyk}
\address{Institute of Mathematics, Polish Academy of Sciences, ul. \'Sniadeckich 8, 00-656 Warszawa, Poland}
\email{b.jakubczyk@impan.pl}

\author{W. Kry\'nski}
\address{Institute of Mathematics, Polish Academy of Sciences, ul. \'Sniadeckich 8, 00-656 Warszawa, Poland}
\email{krynski@impan.pl}
\thanks{W.K. is partially supported by the grant 2019/34/E/ST1/00188 from the National Science Centre, Poland.}

\begin{document}

\begin{abstract}
We study the geometry of dynamic pairs $(X,\VV)$ on a manifold $M$, where $X$ is a vector field and $\VV$ is a distribution on $M$, both satisfying a regularity condition. Special cases are pairs defined by systems of second order ODEs, geodesic sprays in Riemannian, Finslerian and Lagranian geometries, semi-Hamiltonian systems and control-affine systems. Analogs of conjugate points from the calculus of variations are defined for the pair $(X,\VV)$. The main results give estimates for the position of conjugate points in terms of a curvature operator, analogously to the Cartan--Hadamard and Bonet--Myers theorems. Contrary to classical cases, no metric is given a priori, the distribution $\VV$ may be nonintegrable and the curvature operator is defined in terms of $(X,\VV)$.
\end{abstract}

\maketitle

\section{Introduction}
In the present paper we consider \emph{dynamic pairs}  $(X, \VV)$, where $X$ is a vector field on a manifold $M$ and $\VV$ is a distribution on $M$ (a field of tangent $m$-planes), both satisfying a natural regularity condition. Such pairs are naturally defined by geodesic sprays in (pseudo)-Riemannian and Finsler geometry, Lagrangian and Hamiltonian mechanics, in spray geometry, and in geometry of systems of second order ODEs. More generally, they appear in geometric descriptions of control-affine systems. General theory of dynamic pairs was developed in \cite{JK} where a number of geometric objects invariantly assigned to dynamic pairs were introduced. These include a curvature operator, a covariant derivative and invariant metrics. 

Our goal here is to apply these notions in order to study the infinitesimal behaviour of trajectories of $X$, in vicinity of a given one, relative to the distribution $\VV$. To this aim we extend the classical notion of conjugate points from the variational calculus to our setting. Here a basic role is played by the aforementioned curvature operator, denoted $K$. This operator is analogous to the curvature operator which appears in the classical Jacobi equation, often called Jacobi endomorphism. Fiberwise, it is a endomorphism $K_x\colon\VV(x)\to\VV(x)$. All data are assumed to be $C^\infty$-smooth.

Our formalism includes, as particular cases, similar questions deeply studied in Riemann and Finsler geometry and in the geometry of Hamiltonian vector fields on symplectic manifolds in presence of Lagrangian foliations with particular attention put on applications in sub-Riemannian geometry (see \cite{Agr,AB,AgrGam1,AgrGam2, AgrZel,BR,CF}). The theory of the Jacobi vector fields is also well established in the context of systems of second-order differential equations (see \cite{CGM,CM,CMS,HM,JP,Sabau}).  In all the mentioned works distribution $\VV$ is integrable. Furthermore,  there is a canonical metric in the Riemannian, Finslerian and Hamiltonian settings and the operator $K$ is symmetric with respect to the metric. In our case there is no canonical metric and the curvature $K$ is not necessarily symmetric. This makes the problem harder to deal.
%One of the objectives of the paper is to develop a technique which allows to treat non-integrable $\VV$ as well as nonsymmetric cases.

Our motivating example comes from Geometric Control Theory. The most often considered class of non-linear systems are control-affine systems of the form
\begin{equation}\label{ca-system1}
\dot x=X(x)+\sum_{i=0}^mu_iY_i(x),
\end{equation}
where $X,Y_1,\dots,Y_m$ are vector fields on a manifold $M$, $x\in M$ is the state of the system, and  $u_1,\dots,u_m$ are components of the control. Such a system defines a dynamic pair $(X,\VV)$ where
\[
\VV=\spn\{Y_1,\dots,Y_m\}.
\]
We shall assume that $X(x)\not=0$ and $\rk\spn\{Y_1,\dots,Y_m,[X,Y_1],\dots,[X,Y_m]\}=2m$, where $[\cdot,\cdot]$ denotes Lie bracket. These assumptions are often met in applications. Vector field $X$, usually called the \emph{drift} of the system, corresponds to the free motion of the system with no controls. In this regard, our results can be interpreted in terms of stability of the flow generated by $X$ with respect to perturbations tangent to $\VV$. In Section \ref{s-conjugate} we introduce the concept of conjugate points (times) and provide estimates for the first conjugate time. 
As a by-product of our results we also obtain new theorems on the location of conjugate points for general second-order ODEs.  In this context our results strengthen a result in  \cite{HM} concerning a specific case where the curvature has a one-dimensional invariant eigenspace.

In the next Section \ref{s-Dynamic pairs} we recall definitions and basic facts concerning the dynamic pairs. Our exposition is based on \cite{JK}. New results are contained in Sections \ref{s-conjugate} and \ref{s-semi-Hamiltonian}. In particular, in Section \ref{ss-Results} the notion of conjugate points is given and main Theorems \ref{thm1} and \ref{thm2} are stated. Further, we consider Jacobi vector fields in Section \ref{ss-Jacobi fields}. They are used in Section \ref{ss-Proofs} in proofs of the main results. Section \ref{s-semi-Hamiltonian} is devoted to (semi)-Hamiltonian version of dynamic pairs where, additionally to $X$ and $\VV$, a skew-symmetric 2-form $\sigma$ on $M$ is given (possibly, non-closed). Finally, Section \ref{s-Examples} is devoted to examples. In particular we consider mechanical control systems (see \cite{BL,RR}), and a special class of vector second order ODEs related to the so-called dancing construction (studied earlier in \cite{BLN,D,KM}). 

\section{Dynamic pairs, normal frames and curvature}\label{s-Dynamic pairs}
In this section we mainly recall basic notions and constructions introduced before in \cite{JK} and we formulate a number of results used later in the paper. Note that it is shown in \cite{JK} that in the case of the second-order systems of ordinary differential equations our constructions reinterpret the standard approach that goes back to \cite{Cartan,Chern, Kosambi}. A similar approach in this context is also used in \cite{CM,CMS,HM,JP,Sabau}. 
%An important new tool introduced in Section \ref{ss-Invariant metric} is a class of invariant metrics.

\subsection{Dynamic pairs}\label{ss-Dynamic pairs def}
Let $M$ be a smooth differentiable manifold. A \emph{dynamic pair} on $M$ is a pair $(X,\VV)$, where $X$ is a vector field on $M$ and $\VV\subset TM$ is smooth distribution on $M$ which is a field of tangent planes of dimension $m=\rk \VV$. The vector fields $V\in\VV$, together with the Lie brackets $[X,V]$ span a new distribution $\DD\subset TM$ denoted
\[
\DD=\VV+[X,\VV].
\]
Above and later on we often write $V\in\VV$, instead of $V\in\Sec(\VV)$, meaning that $V$ is a section of $\VV$, usually
a local section.

\begin{definition}\label{d2.0}
A dynamic pair $(X,\VV)$ is called \emph{regular} on $U\subset M$ if it satisfies three conditions:
\[
X(x)\neq0, \quad x\in U, \eqno{(R1)}
\]
\[
\rk \DD=2m \qquad\qquad \eqno{(R2)}
\]
(i.e. the rank of $\DD$ is maximal possible), and
\[
[X,\DD]\subset\DD, \eqno{(I)}
\]
i.e. $\DD$ is invariant under the local flow of $X$. Conditions (R1) and (R2) will be called \emph{regularity conditions}, while (I) will be called \emph{invariance condition}.
\end{definition}

\begin{remark}
We will usually identify $M$ with $U$, taking it sufficiently small so that (R1), (R2) and (I) are satisfied on $M$. A sufficiently interesting case is that of $D=TM$. In this case $\dim M=2m$ and the invariance condition (I) is automatically satisfied. Another natural case is $\dim M= 2m+1$ and $TM=\DD\oplus\spn\{X\}$. In this case most results of the paper also hold under a weaker condition: $[X,\DD]\subset\DD\mod X$ which only requires that all computations are performed modulo $X$. However, for simplicity, we assume that (I) is always satisfied. Condition (R1) is essential in Proposition \ref{p2.2}, where it implies the non-singularity of equation \eqref{eq5}.
\end{remark}

\subsection{Curvature operator}\label{ss-Curvature}
If $V_1,\ldots,V_m$ span $\VV$ on $M$ or an open subset $U\subset M$ (resp. along a trajectory $c$ of $X$) then the tuple 
\[
\hat V=(V_1,\ldots,V_m),
\]
treated as a row vector, will be called a \emph{frame} of $\VV$ on $U$ (resp. along $c$).

Note that condition (R2) implies that the vector fields $V_j$ and $[X,V_j]$, $0\le j\le m$, are pointwise linearly independent and they span $D$.
Moreover, it follows from (I) that $[X,[X,V_j]]$ are also in $D$. Thus 
\[
[X,[X,V_j]]=\sum_i(H_0)^i_jV_i+\sum_i(H_1)^i_j[X,V_i],
\]
where the functions $(H_0)^i_j$, $(H_1)^i_j$ form two $m\times m$ matrices $H_0$ and $H_1$. In the matrix notation we have
\begin{equation}\label{eq1}
[X,[X,\hat V]]=\hat V H_0+[X,\hat V]H_1.
\end{equation}

\begin{definition}\label{d2.1}
The curvature matrix in the frame $\hat V$, of the regular pair $(X,\VV)$, is the matrix
\begin{equation}\label{eq-curvature}
K=-H_0+\frac{1}{2}X(H_1)-\frac{1}{4}H_1^2,
\end{equation}
where $X(H_1)$ is the matrix of derivatives along $X$ of the coefficients of $H_1$.
\end{definition}

The following result proved in \cite[Formula (11) and Proposition 2.5]{JK} implies that the matrix $K$ defines, at each $x\in M$, a linear operator $K_x:\VV(x)\to\VV(x)$ which will be called \emph{curvature operator} of $(X,\VV)$.

\begin{proposition}\label{p2.1}
When the local frame $\hat V$ is changed for $\hat V'=\hat V G$, with an  $m\times m$ invertible matrix $G$ of smooth functions, the curvature matrix is
transformed to
\begin{equation}\label{eq-curvature-transformed}
K'=G^{-1}KG.
\end{equation}
\end{proposition}
The curvature $K$ can be defined in a more conceptual way. For this we need the following notion.

\begin{definition}\label{d2.2}
A local section $V$ of $\VV$  is called \emph{normal} if
\[
[X,[X,V]] \in \VV.
\]
A local frame $\hat V=(V_1,\ldots,V_m)$ of $\VV$ on $U\subset M$ is called \emph{normal frame} of $\VV$ if all $V_i$ are normal sections of $\VV$. We also call $V_1,\ldots,V_m$ \emph{normal
generators} of $\VV$.
\end{definition}

Normal frames in $\VV$ exist (see \cite[Proposition 2.2]{JK}) and have the following properties.

\begin{proposition}\label{p2.2}
For a regular dynamic pair $(X,\VV)$ on $M$ the following holds. 

(i) Given $x_0\in M$ and a frame $F$ in $\VV(x_0)$, then along the (non-closed) trajectory $c$ of $X$ starting from $x_0$ there is a unique normal frame $\hat V$ of $\VV$ which coincides with $F$ at $x_0$. In particular, there is a normal frame in neighbourhood of $x_0$. 

(ii) If $\hat V$ is a frame of $\VV$ on an open $U\subset M$ and $H_1$ is defined via \eqref{eq1}, then $\hat W=\hat V G$ is a normal frame of $\VV$ on
$U$ if and only if
\begin{equation}\label{eq5}
X(G)=-\frac{1}{2}H_1G.
\end{equation}
In particular, if $\hat V$ is normal then $\hat W=\hat V G$ is normal if and only if $X(G)=0$ and a vector field $Y=\sum f^iV_i$ is normal if and only if $X(f^i)=0$.
\end{proposition}

\begin{remark}
Equation \eqref{eq5} implies that the function $d(t)=\det G(x(t))$ satisfies the differential equation $\dot d=ad$ along any trajectory $x(t)$ of $X$, where $a(t)=-\frac{1}{2}\tr H_1(x(t))$. Assuming $d(0)\not=0$ we have $d(t)\ne0$, thus \eqref{eq5} has a global, non-degenerate solution $G$ along the whole trajectory $x(t)$.
\end{remark}

Now, if (R1), (R2), (I) hold and $\hat V=(V_1,\ldots,V_m)$ is a normal frame in $\VV$, there are functions $K_i^j$ such that
\begin{equation}\label{eq13}
[X,[X,V_i]]+\sum_jK^j_iV_j=0,
\end{equation}
or in the matrix form
\begin{equation}\label{eq14}
[X,[X,\hat V]]+\hat V K=0.
\end{equation}
This means that $H_0=-K$ and $H_1=0$ in equation \eqref{eq1}. Thus $K=(K_i^j)$ coincides with the curvature matrix defined in Definition \ref{d2.1}. The fact that the curvature matrix $K$ is defined in a normal frame will have a special meaning for us. Thus we introduce

\begin{definition}\label{d2.3}
Let $c:[0,T]\to M$ be an integral curve of $X$. The curvature matrix
$$
K(t)=(K_i^j(c(t))), \quad t\in[0,T],
$$
expressed in a normal frame along $c$ (equivalently, defined by equation \eqref{eq13} along $c$) will be called \emph{normal curvature matrix} along $c$
of the regular pair $(X,\VV)$.
\end{definition}

\subsection{Systems of second order ODEs}\label{ss-ODEs}
As mentioned before, basic examples of regular dynamic pairs are provided by systems of second order differential equations $\ddot x=F(t,x,\dot x)$, $t\in \R$, $x\in \R^m$. Geometric invariants of such systems were deeply studied by many authors beginning with the works of Kosambi \cite{Kosambi}, Cartan \cite{Cartan}, and Chern \cite{Chern}. Our general constructions presented in the next sections coincide, when specified to this case, with works of \cite{CM,JP} (see also \cite{HM} and references therein). An invariant characterization of general dynamic pairs which correspond to systems of ODEs is given in \cite[Theorem 4.5]{JK}.

For systems of second order ODEs we take $M=J^1(\R,\R^m)\simeq\R\times\R^m\times\R^m$ which is the space of 1-jets of mappings $\R\to\R^m$ with natural coordinates $t,x^i,y^i$ on $M$ ($y^i$ replacing $\dot x^i$). We have a natural vector field on $M$ which is the total derivative
\[
X=\partial_t+\sum_{i=1}^m y^i\partial_{x^i}+\sum_{i=1}^m F^i\partial_{y^i}.
\]
As the distribution $\VV$ we take $\spn\,\{\partial_{y^1},\ldots,\partial_{y^m}\}$ (the tangent distribution to the foliation of fibers of the bundle of
1-jets of curves in $\R^m$). It is then easy to check that, for $V_j=\partial_{y^j}$, we have $[X,V_j] = -\partial_{x^j} - \sum_iF^i_{y^j}\partial_{y^i}$
and
\[
\DD=\spn\,\{\partial_{x^1},\ldots,\partial_{x^m}, \partial_{y^1},\ldots,\partial_{y^m}\}
\]
Thus $(X,\VV)$ is a regular dynamic pair. Moreover, we have
\[
[X,[X,V_j]] = \sum_iF^i_{y^j}\partial_{x^i} + \sum_{i,k}(-F^i_{ty^j}-y^kF^i_{x^ky^j}-F^kF^i_{y^ky^j}+F^i_{x^j}+F^k_{y^j}F^i_{y^k})\partial_{y^i}.
\]
Taking into account that the first sum and the last term in the second sum give together $-\sum_iF^i_{y^j}[X,V_i]$ we see that the matrices in
\eqref{eq1} are equal to
\[
H_0=(-F^i_{ty^j}-y^kF^i_{x^ky^j}-F^kF^i_{y^ky^j}+F^i_{x^j}), \quad H_1=(-F^i_{y^j}).
\]
Thus the general formula \eqref{eq-curvature} for curvature gives the matrix 
\begin{equation}\label{curvatureODE'}
K^i_j = -F^i_{x^j} - \frac{1}{4}F^i_{y^k}F^k_{y^j} + \frac{1}{2}F^kF^i_{y^ky^j} + \frac{1}{2}y^kF^i_{x^ky^j} + \frac{1}{2}F^i_{ty^j},
\end{equation}
or in more compact form
\begin{equation}\label{curvatureODE}
K=-F_x+\frac{1}{2}X (F_{\dot x})-\frac{1}{4}(F_{\dot x})^2.
\end{equation}
Formula \eqref{curvatureODE'} appeared for the first time in the aforementioned works \cite{Cartan, Chern,Kosambi} and, in recent literature, it is often referred to as the Jacobi curvature or Jacobi endomorphism (cf. \cite{CM}). Note that the trace-free part of $K$ is referred to as the torsion in \cite{Grossman}.

\begin{remark} 
Notice, that the case of autonomous equations $\ddot x=F(x,\dot x)$ can also be analysed on the simpler manifold $M=\R^m\times\R^m$, in which case we remove $\partial_t$ from the total derivative. The remaining ingredients define a regular dynamic pair with the distribution $\DD$ being the whole tangent bundle $TM$ and the curvature given by \eqref{curvatureODE'}, with $F_{ty}=0$.
\end{remark}

\subsection{Invariant splitting and connection}\label{ss-Splitting}
Proposition \ref{p2.2} says that if both $\hat V$ and $\hat W=\hat V G$ are normal frames of $\VV$ then $X(G) = 0$. Hence $[X,\hat W] =[X,\hat V]G$. Consequently, given a normal frame $\hat V=(V_1,\dots,V_m)$, the distribution
\[
\HH=\spn\{[X,V_1],\ldots,[X,V_m]\}=\spn\{[X,\hat V]\}
\]
does not depend on the choice of such frame.  Clearly $\HH\subset\DD$, which follows from $\VV\subset \DD$ and $[X,\DD]\subset \DD$.

The distribution $\HH$ will be called \emph{horizontal} and vector fields in $\HH$ will also be called \emph{horizontal}. By Proposition \ref{p2.2} a vector field $Y=\sum f^iV_i$, expressed in a normal frame $(V_1,\dots,V_m)$, is normal iff $X(f^i)=0$ and then $[X,Y]=\sum f^i[X,V_i]$ is horizontal. Conversely, if $\sum f^i[X,V_i]$ is
horizontal then $X(f^i)=0$.

If $\hat V=(V_1,\dots,V_m)$ is a general frame of $\VV$, then $\HH$ can be computed as follows (see \cite[Proposition 2.5]{JK}).

\begin{proposition}\label{p-horizontal}
If $[X,[X,\hat V]]=[X,\hat V]H_1+\hat V H_0$, for some $m\times m$ matrices $H_0$, $H_1$,  then
\[
\HH=\spn\left\{\ad_X\hat V-\frac{1}{2}\hat V H_1\right\}.
\]
\end{proposition}

Condition (R2) implies that the distribution $\DD$ can be written as
\[
\DD=\VV\oplus\HH
\]
which we call the \emph{canonical split} of $\DD$. This relation defines canonical projections
\[
\pi_\VV\colon \DD\to\VV\quad \mathrm{and}\quad \pi_\HH\colon \DD\to\HH.
\]
\begin{definition}\label{d4.3}
A \emph{covariant derivative} along $X$ of vector fields in $\VV$ is the operator $D_X\colon\Sec(\VV)\to\Sec(\VV)$ given by
\[
D_XV:=\pi_\VV[X,V], \quad V\in\VV.
\]
More generally, $D_X$ extends to $D_X\colon\Sec(\DD)\to\Sec(\DD)$ by formula
\[
D_XY:= \pi_\VV[X,\pi_\VV Y] + \pi_\HH [X,\pi_\HH Y].
\]
\end{definition}

Note that $D_X$ is also well defined as an operator $\Sec_c(\VV)\to\Sec_c(\VV)$, or $\Sec_c(\DD)\to\Sec_c(\DD)$, where $\Sec_c$ denotes sections along a given trajectory $c\colon[0,T]\to M$ of $X$. It satisfies the Leibniz rule
\[
D_X(fY)=X(f)V+fD_XY.
\]
and the definition implies the following characterization.
\begin{proposition}\label{p4.4} For vector fields $V\in \VV$ and $H\in\HH$ we have:
\[
D_XV=0 \ \Longleftrightarrow\ V\ \it{is\ normal}, \qquad\qquad\qquad\ \
\]
\[
\qquad\qquad D_XH=0 \ \Longleftrightarrow\ H=[X,W],\ \it{where}\ W\in\VV\  \it{is\ normal}.
\]
\end{proposition}

\begin{remark}
The canonical splitting in the context of second order systems of ODEs and the corresponding covariant differential $D_X$ coincides with the analogous notions used in \cite{CM,JP}. Counterparts of these notions for higher order ODEs were introduced in \cite{JK}.
\end{remark}

\begin{remark}\label{r3}
The operators $A\colon\VV\to\HH$ and $B\colon\HH\to\VV$ defined by
\[
 A=\pi_\HH\circ\ad_X, \quad B=\pi_\VV\circ\ad_X
\]
are vector bundle morphisms, canonically defined by $(X,\VV)$. This is a consequence of the fact that, any time we compose the operator $\ad_X$ defined on sections $Sec(\Delta)$ of a distribution $\Delta\subset TM$ with a projection $\pi$ along $\Delta$, we obtain a vector bundle morphism because $\pi(\ad_X(fY))=\pi(X(f)Y+f[X,Y])=f\pi([X,Y])$. In addition, it follows from (R2) that $A$ is a local isomorphism of bundles, given on a normal frame by $AV_i=[X,V_i]$. The operators $A$ and $B$ have canonical extensions to $\DD$ defined by $\hat A=\pi_\HH\circ\ad_X\circ\pi_\VV$ and $\hat B=\pi_\VV\circ\ad_X\circ\pi_\HH$.
Using them we get the decomposition $\ad_X=D_X+\hat A+\hat B$.

\end{remark}
\begin{remark}
The curvature operator $K$ of a regular dynamic pair $(X,\VV)$ can be equivalently defined as $K = - B\circ A$ which follows from formula \eqref{eq13}.
\end{remark}

\subsection{Invariant metrics}\label{ss-Invariant metric}
Any regular dynamic pair $(X,\VV)$ admits, along any trajectory of $X$ (or locally), metrics which are invariant under covariant differentiation along $X$. We begin with defining such a metric on the distribution $\VV$. Consider a normal frame $\hat V=(V_1,\dots,V_m)$ in $\VV$, along a trajectory $c\colon[0,T]\to M$ of $X$. Let $g$ be the positive definite metric in $\VV$ such that $V_1,\dots,V_m$ are orthonormal. Equivalently, if $Y=\sum \phi^iV_i$, $Y_2=\sum \psi^iV_i$ are vector fields in
$\VV$, along $c$, the metric $g$ is given by
\[
g(Y,Z)=\sum_i\phi^i\psi^i.
\]
By adding minus signs in the above sum we can define a metric with any prescribed signature.
It follows from Proposition \ref{p4.4} that $D_XV_i=0$, thus $D_XY=\sum X(\phi^i)V_i$ and $D_XZ=\sum X(\psi^i)V_i$ and we see that the metric is invariant under the covariant derivative $D_X$, i.e.,
\begin{equation}\label{eq-der}
X(g(Y,Z))=g(D_XY,Z)+g(Y,D_XZ).
\end{equation}

Any metric on $\VV$ satisfying \eqref{eq-der} will be called \emph{invariant metric}. The above metric on $\VV$ can be extended to $\DD=\VV\oplus\HH$.
Namely, given a normal frame $\hat V=(V_1,\dots,V_m)$ in $\VV$, we define $g$ by declaring that the vector fields $V_1,\dots,V_m$, $[X,V_1],\dots,[X,V_m]$ are orthonormal. Then, for arbitrary vector fields in $\DD$,
\[
Y=\sum \phi^iV_i+\sum \widetilde \phi^i[X,V_i], \quad Z=\sum \psi^iV_i+\sum \widetilde \psi^i[X,V_i]
\]
(all defined along $c$ or locally) we have
\[
g(Y,Z)=\sum_i(\phi^i\psi^i+\widetilde \phi^i\widetilde \psi^i).
\]
By Proposition 2.12 we also have $D_X[X,V_i]=0$, thus the invariance formula \eqref{eq-der} holds for all $Y,Z\in \DD$. Note that, since the normal frame along $c$ is uniquely defined by its initial value at $c(0)$, by Proposition 2.6, the invariant metric along $c$ is unique up to a linear transformation in $\DD(c(0))$. Thus we have proved

\begin{proposition}\label{p-x}
Suppose we are given a regular dynamic pair $(X,\VV)$, an integral curve $c:[0,T]\to M$ of $X$, and a frame $F$ in $\VV(c(0))$. Then the following holds.

(i) There is a unique invariant positive definite (respectively, of prescribed signature) metric $g$ on $\DD$ along $c$ such that $F$ is orthonormal with respect to $g$ restricted to $\VV(c(0))$. The normal frame which coincides with $F$ at $c(0)$ is orthonormal with respect to $g$ along $c$.

(ii) $\VV$ and $\HH$ are orthogonal with respect to $g$ along $c$. 

(iii) $g(AV,AW)=g(V,W)$ for all $V,W\in \VV$ along $c$. 
\end{proposition}

Given a metric $g$ on $\VV$, possibly not positive definite, we can define a \emph{directional curvature} of a regular dynamic pair $(X,\VV)$, using its curvature operator $K$.
\begin{definition}\label{d-directional curvature}
The \emph{curvature  in direction} $v\in \VV$ of the triple $(X,\VV,g)$ is
\[
k(v)=\frac{g(Kv,v)}{g(v,v)}, \quad \mathrm{for}\quad v\in\VV\quad\mathrm{such\ that} \quad g(v,v)\ne0.
\]
\end {definition}

\section{Conjugate points and Jacobi vector fields}\label{s-conjugate}
In this section we introduce a definition of conjugate points for dynamic pairs and state our main results on existence and location of conjugate points in terms of the curvature operator of the dynamic pair. For proving the results we introduce Jacobi vector fields generalized to the current context and establish correspondence between conjugate points and properties of Jacobi fields (Section \ref{ss-Jacobi fields}). The proofs follow in Section \ref{ss-Proofs}.

\subsection{Conjugate points}\label{ss-Results}
Let $(X,\VV)$ be a regular dynamic pair. Consider an integral curve $c\colon[0,T]\to M$ of $X$, $c(0)=x_0$, which is not a closed orbit. Regularity condition (R2) implies that the space $\VV(x_0)$ transported by the flow of $X$ along the trajectory $c(t)$ is transversal to the space $\VV(c(t))$ for $t$ small enough. Points $c(t)$ where this is not true will be called conjugate to $c(0)$, analogously as in second order ODEs and Calculus of Variations.

\begin{definition}\label{d3.1}
Given $c$, a point $0<t^*\le T$ is called \emph{conjugate time} and the corresponding point $c(t^*)$ \emph{conjugate point} to $x_0$ on $c$ if
$$
\exp(t^*X)_*\VV(x_0)\cap \VV(c(t^*))\not=\{0\}.
$$
The dimension of the above intersection is called \emph{the multiplicity} of $t^*$ or $c(t^*)$.
\end{definition}

The intuition behind the definition is  as follows. Consider a smooth family  $c_s\colon[0,T]\to M$ of integral curves of $X$, $s\in (-\epsilon,\epsilon)$, uniquely defined by the curve of initial points $\gamma(s)=c_s(0)$. Then $t^*>0$ is conjugate time if there exists a curve of initial points $\gamma(s)$ which for all $s$ is tangent to $\VV$ such that the curve of endpoints $s\to c_s(t^*)$ is tangent to the subspace $\VV(c_0(t^*))$, i.e.  $\frac{d}{ds}c_s|_{s=0}(t^*)\in\VV(c_0(t^*))$. Informally, there exists a vector $v\in\VV(c(0)$ such that an infinitesimally small perturbation of the initial condition $c(0)$ in the direction of $v$ causes the perturbation of the trajectory $c(t)$ of $X$ which, after time $t^*$, is is infinitesimally small in a direction $w\in \VV(c(t^*))$'. 

If $\VV$ is the control distribution $\spn\{Y_1,\dots,Y_m)$, defined by a control affine system \eqref{ca-system1}, then this means that there is a control value $u$ such that if a small control impulse in direction $u$ is applied at time $t=0$ then the perturbed trajectory arrives at time $t^*$ to a point infinitesimally close to the reference trajectory in a direction which belongs to $\VV(c(t^*))$, thus another impulse of control can ``remove'' the effect of the impulse at $t=0$, placing the perturbed trajectory back to $c(t)$.  

\begin{remark}
If distribution $\VV$ is integrable then one can consider locally defined quotient manifold $N=M/\VV$ that comes with a natural projection $\pi\colon M\to N$. Then any integral curve $t\mapsto c(t)$ of $X$  defines a curve in $N$ via projection  $\tilde c(t)=\pi(c(t))$. If a family of curves $c_s(t)$, $s\in (-\epsilon,\epsilon)$, satisfies $\frac{d}{ds}c_s(0)\in\VV$ for all $s$ then $\tilde c_s(0)=\tilde c_0(0)$ for all $s$ and the curve $s\to c_s(0)$ defines a point in $N$. Hence we arrive to the usual definition of conjugate points in this case. 
\end{remark}

Now we can formulate two results that generalize the classical Cartan--Hadamard and Bonnet--Myers theorems on conjugate points on Riemannian and Finslerian manifolds.

\begin{theorem}\label{thm1}
Given a regular dynamic pair $(\VV,X)$ and an invariant metric $g$, if there exists $\lambda\in\R$ such that one of the following equivalent estimates holds
\begin{equation}\label{estimate-lambda1}
g(Kv,v) \le \lambda g(v,v), \quad \mathrm{for}\ \  v\in \VV(c(t)),\ \ t\in [0,T],
\end{equation}
or, using directional curvature,
\begin{equation}\label{eq-dir-curv-estimate}
k(v)\leq \lambda, \quad v\in \VV(c(t))\setminus\{0\}, \ \ t\in [0,T], 
\end{equation}
then there are no conjugate times in $(0,T]$, if $\lambda\le 0$, and there are no conjugate times in the interval $(0,t_c)$, $t_c=\min\{T,\frac{\pi}{\sqrt{\lambda}}\}$, if $\lambda>0$.
\end{theorem}

\begin{remark}
Note that above estimates can be replaced by the equivalent estimate
\begin{equation}\label{estimate-lambda2}
\lambda^K(t) \le \lambda, \quad t\in [0,T],
\end{equation}
where $\lambda^K(t)$ denotes the maximal eigenvalue of the symmetric (with respect to the scalar product $g(t)$) part of $K(t)$.
\end{remark}

\begin{remark}
If the curvature $K(t)$ is diagonalizable and the corresponding eigenspaces are parallel with respect to $D_X$ then the above result (in the case of ODEs) reduces to a result in \cite{HM}. 
\end{remark}

To state the second result consider a non-closed integral curve $c:[0,T)\to M$ of $X$, with $T\in (0,\infty]$, and let $(X,\VV)$ be a regular dynamic pair defined in a neighbourhood of $c$. 

\begin{theorem}\label{thm2}
If, along $c$, the curvature $K$ is symmetric with respect to some, parallel with respect to $D_X$, positive definite metric $g$ and
\[
\tr K\ge \kappa>0\ \ \text{along}\ c
\]
then there is a conjugate point in the interval $(0,T^*]$, $T^*=\pi\sqrt{\frac{m}{\kappa}}$, provided that $T^*<T$.
\end{theorem}

\begin{remark}
The symmetry of $K$ in the above theorem is crucial, as demonstrated in Example \ref{s-Example} at the end of the paper.
\end{remark} 

The following result does use any metric and, thus, it does not require symmetry of $K$.

\begin{theorem}\label{thm3}
(i) If the curvature matrix $K_t=K(c(t))$, $t\in [0,T]$, has an invariant one-dimensional subspace
which is parallel with respect to $D_X$ along $c$ and the corresponding eigenvalue satisfies 
\begin{equation}\label{LA}
\lambda(t)\ge \kappa>0, \ \ t\in(0,T], 
\end{equation}
(equivalently, the normal curvature matrix $K_t$ on $c$ has a one-dimensional diagonal block with the eigenvalue satisfying \eqref{LA}) 
then there is a conjugate time in $(0,\pi/\sqrt{\kappa}\,]$, provided that $\pi/\sqrt{\kappa}\le T$. 

(ii) Moreover, if there are $r$ such subspaces with the eigenvalues $\lambda_1(t),\dots,\lambda_r(t)$ then the set of conjugate points along $c$ of the pair $(X,\VV)$ includes all zeros on $(0,T]$ of solutions of the equations
\begin{equation}\label{SJ}
\ddot y_i=-\lambda_i(t)y_i, 
\end{equation}
$i=1,\dots,r$, with initial conditions $y_i(0)=0$, $\dot y_i(0)\not=0$ (if $t^*$ is a common zero of $k$ of these equations then it is conjugate point of multiplicity $k$ or larger) .
\end{theorem}

\begin{remark}
The assumption on eigenspaces being parallel with respect to $D_X$ is equivalent to the fact that the normal curvature matrix is diagonalizable. In Example \ref{ex-SODE} we will give a geometric characterization of 2-dimensional systems admitting such decomposition. 
\end{remark}

\subsection{Jacobi vector fields}\label{ss-Jacobi fields}
Consider a regular dynamic pair $(X,\VV)$ on $M$. In order to prove main results on conjugate points we extend the classical notion of the Jacobi fields to our context. It will have two variants. 

 \begin{definition}\label{d6.1}
A vector field $V\in\VV$ is a \emph{vertical Jacobi field} of $(X,\VV)$ if
\begin{equation}\label{eq-JE}
D_XD_XV+KV=0.
\end{equation}
A vector field $J$ in $\DD$ is \emph{a Jacobi field} of $(X,\VV)$ if
\begin{equation}\label{eq-J1}
[X,J]=0.
\end{equation}
Both definitions make sense on an open subset of $M$, as well as on a single integral curve $c:[0,T]\to M$ of $X$. In the latter case we replace
condition $[X,J]=0$ by
\begin{equation}\label{eq-J1'}
\exp(tX)_*J_0=J_t, \quad t\in [0,T],
\end{equation}
where $J_t$ is the vector field $J$ at $c(t)$. Equation \eqref{eq-JE} is called \emph{Jacobi equation}.
 \end{definition}

Condition \eqref{eq-J1} implies that the flow $\exp(tX)$ of $X$ transports trajectories of the Jacobi field $J$ onto its trajectories, preserving time
parametrization. In particular, given an integral curve $\gamma(s)$ of $J$,  the curve $\gamma_t(s)=\exp(tX)(\gamma(s))$ is again a trajectory of $J$. In
studying conjugate points of $(X,\VV)$ we are particularly interested in curves $\gamma$ passing through $x_0=\gamma(0)$ which are tangent to $\VV$ at
$x_0$ and transporting them by the flow $\exp(tX)$.

Jacobi fields and vertical Jacobi fields are in one-to-one correspondence. To show this we use the vector bundle isomorphism $A:\VV\to\HH$, $AV=\pi_\HH[X,V]$ (Remark \ref{r3}).

 \begin{proposition}\label{p6.2}
If  $V\in\VV$ is a vertical Jacobi field of $(X,\VV)$ then the vector field
\begin{equation}\label{eq-J-up}
J_V=-D_XV+AV
\end{equation}
is a Jacobi field of $(X,\VV)$.
Vice versa, if  $J$ is a Jacobi field of $(X,\VV)$ then the vector field
\begin{equation}\label{eq-J-down}
V_J=A^{-1}\pi_\HH J
\end{equation}
is a vertical Jacobi field of $(X,\VV)$.
 \end{proposition}

\paragraph{\bf Proof.}
Let $V$ be a vertical Jacobi field. Locally, we may write $V=\sum f^iV_i$, where $V_1,\dots,V_m$ are normal generators of $\VV$. Then $D_XV_i=0$ and $D_XV=\sum X(f^i)V_i$. From the definition of $A$ and the fact that $[X,V_i]\in\HH$, by the definition of $\HH$, we have $AV=\sum f^i[X,V_i]$. Thus
\begin{align*}
[X,J_V]&=-[X,\sum X(f^i)V_i]+[X,\sum f^i[X,V_i]]\\
&=-\sum X^2(f^i)V_i +\sum f^i[X,[X,V_i]]=-D_XD_XV-KV=0,
\end{align*}
since $D_XD_X(\sum f^iV_i)=\sum X^2(f^i)V_i$ and $K(\sum f^iV_i)=\sum_{i,j}f^iK_i^jV_j=-\sum_if^i[X,[X,V_i]]$.

In order to prove the second assertion, let $V_1,\dots,V_m$ be a local normal frame of $\VV$. Given a Jacobi field $J$, we may write
\[
J=\sum_ig^iV_i+\sum_if^i[X,V_i].
\]
Then
\[
V_J=A^{-1}\pi_\HH J=A^{-1}\left(\sum f^i[X,V_i]\right)=\sum f^iV_i.
\]
We will show that $V_J=\sum f^iV_i$ is a vertical Jacobi field. The assumption $[X,J]=0$ gives
\[
0=\sum (X(g^i)V_i+f^i[X,[X,V_i]])+\sum (g^i+X(f^i))[X,V_i]
\]
and, taking into account $[X,[X,V_i]]=\sum_j-K_i^jV_j$, we see that $[X,J]=0$ is equivalent to
\begin{equation}
g^i+X(f^i)=0,\qquad X(g^i)-\sum_jK^i_jf^j=0.
\end{equation}
 Thus $X^2(f^i)=-\sum_jK^i_jf^j$ and we find that
\[
D_X^2\left(\sum f^iV_i\right)=\sum X^2(f^i)V_i=-\sum_{i,j}f^jK^i_jV_i=-\sum f^j KV_j = -K\left(\sum f^jV_j\right),
\]
i.e., $D_X^2V_J=-KV_J$. This means that $V_J$ is a vertical Jacobi field.
\QED

Clearly, existence of a vector field $J$ in a neighbourhood of $c$ which satisfies condition \eqref{eq-J1} implies existence of a vector field $J_t$ along $c$ satisfying \eqref{eq-J1'}, and vice versa. According to its definition a point $t^*\in (0,T]$ is conjugate time if
there exists a nonzero $v=J_0\in \VV(x_0)$ such that $\exp(t^*X)_*J_0\in \VV(c(t^*))$ or, equivalently, there exists a nonzero Jacobi field $J_t$ along $c$ such that
\begin{equation}\label{eq-J2}
J_0\in \VV(c(0)) \quad \mathrm{and} \quad J_{t^*}\in \VV(c(t^*)).
\end{equation}
Using Proposition \ref{p6.2} we see that existence of a Jacobi field $J_t$ along $c$ is equivalent to existence of a vertical Jacobi field $V_t$ and $V_t=A^{-1}\pi_\HH J_t$, in which case the above boundary conditions are equivalent to $V_0=0$ and $V_{t^*}=0$. Additionally, the dimensions of the space of such Jacobi vector fields and the corresponding vertical Jacobi vector fields coincide.
Thus we obtain

\begin{proposition}\label{p6.3}
A point $t^*\in (0,T]$ is a conjugate time along an integral curve $c\colon[0,T]\to M$ of $X$ iff there exists a nontrivial vertical Jacobi field $V\in\VV$ along $c$ such that
\begin{equation}\label{eq-J3}
V(c(0))=0 \quad \mathrm{and} \quad V(c(t^*))=0.
\end{equation}
The multiplicity of a conjugate time equals to the number of linearly independent vertical Jacobi fields satisfying \eqref{eq-J3}.
\end{proposition}

\subsection{Proofs}\label{ss-Proofs}
With all earlier preparations the proofs can be reduced to classical arguments.  As in the theorem, we will use a normal frame in $\VV$ along $c$ and the corresponding invariant metric $g$ on $\VV$ making the frame orthonormal (see the beginning of Section \ref{ss-Invariant metric}). The curvature matrix $K(t)$, called normal, is taken in such frame.
\vskip 2ex

\paragraph{\bf Proof of Theorem \ref{thm1}.}
Consider a vertical Jacobi field $V\in \VV$ along $c:[0,T]\to M$, with $V(c(0))=0$, $(D_XV)(c(0))=v\not=0$. By Proposition \ref{p6.3} it is enough to estimate the first zero
of $V_t=V(c(t))$ in $(0,T]$. This is equivalent to estimating the first zero of the function
\[
\varphi=\sqrt{g(V,V)}
\]
on $c(t)$,  $t\in (0,T]$. From the invariance of the metric \eqref{eq-der} we have
\[
\varphi'=\frac{g(D_XV,V)}{\varphi},
\]
and
\[
\varphi''=\frac{g(D_XD_XV,V)}{\varphi} + \frac{g(D_XV,D_XV)}{\varphi} - \frac{g^2(D_XV,V)}{\varphi^3}.
\]
The Cauchy--Schwartz inequality $g^2(D_XV,V) \le g(D_XV,D_XV)\ g(V,V)$ and $g(V,V)=\varphi^2$ imply that the last two terms in the above expression give a nonnegative number. This and the Jacobi
equation for $V$ imply that
\[
\varphi'' \ge -\frac{g(KV,V)}{\varphi}=-\frac{k(V)g(V,V)}{\varphi}=-k(V)\varphi.
\]
From the assumption $\lambda^K\leq \lambda$, equivalent to the estimate \eqref{eq-dir-curv-estimate}, we finally get
\begin{equation}\label{eq-inequality}
\varphi'' \ge -\lambda \varphi.
\end{equation}
Note that the initial conditions on $V$ imply that $\varphi(0)=0$ and $\varphi'(0)=\sqrt{g(v,v)}>0$ (since, for $\psi=\varphi^2$, we have $\psi'=2\varphi'\varphi=2g(D_XV,V)$, thus   $\psi(0)=\psi'(0)=0$, and $\psi''=g(D^2_XV,V)+g(D_XV,D_XV)$, thus $\psi''(0)=g(v,v)$). 
It follows that $\varphi$ is positive on some interval $(0,t_z)$, where $t_z$ is taken the first zero of $\varphi$ in $(0,T]$, or $t_z=T$ if $\varphi$ does not have such a zero. Now we can write, instead of
\eqref{eq-inequality}, the equality
\begin{equation}\label{eq-equality}
\varphi'' =a(t) \varphi
\end{equation}
on $(0,t_z)$, where $a(t)=\varphi''/\varphi\ge -\lambda$. Applying the classical Sturm comparison theorem (see \cite{DMJ}) we find out that the first zero $t^*$ of $\varphi$ is larger or equal to the first zero of the equation $\psi''=-\lambda\psi$ with $\psi(0)=0$, $\psi'(0)=\varphi'(0)>0$, thus $t^*\ge \pi/\sqrt{\lambda}$, if $\lambda>0$, or $t^*$ does not exist in $(0,T]$.
\QED\vskip 2ex

For the proof of Theorem \ref{thm2} we need a lemma. Let $g$ be a positive definite metric on $\VV$ defined along $c$ and parallel with respect to $D_X$. Given $0<r< T$, we introduce the \emph{index functional} along $c$ defined on sections $V$ of $\VV$ along $\bar c=c\vert_{[0,r]}$ by
\begin{equation}\label{I1}
I(V,V)=\int_0^{r}( g(D_XV,D_XV)-g(KV,V))dt. 
\end{equation}

\begin{lemma}\label{lem-index1}
Let $K$ be symmetric with respect to $g$ on $c$ and let $0<r<T$. Then there are no conjugate points in $(0,r]$ if and only if 
for any section $V$ of $\VV$ along $c$ such that $V(c(0))=0=V(c(r))$ we have
\begin{equation}\label{I2}
I(V,V)\ge 0 
\end{equation}
and equality holds only if $V\equiv0$ on $[0,r]$.
\end{lemma}

The proof will follow from Lemma \ref{lem-index2}
\medskip

\paragraph{\bf Proof of Theorem \ref{thm2}.}
Assume that there are no conjugate points in the interval $(0,T^*]$. We will show that this leads to contradiction with the definition of $T^*$. Consider the interval $I=[0,r]$ where $T^*<r<T$. Then there there are no such points in $(0,r]$ for some $r>T^*$ close enough to $T$. This follows from the fact that $r$ is a conjugate time if and only if the m-tuple $J_1,\dots,J_m$ of Jacobi fields satisfying the initial conditions $J_i(0)=0$ and $D_XJ_i(0)=e_i$, with $e_1,\dots,e_m$ a basis in $\VV(c(0))$, is linearly dependent at $c(r)$. Linear independence of the tuple at $c(T^*)$ implies their independence at $c(r)$ for $r$ close o $T^*$ which implies the above statement.
 
Let $(V_1,\dots,V_m)$ be a normal frame along $c$ defined by an orthonormal with respect to $g$ basis $e_1,\dots,e_m\in \VV(c(0))$ and the initial conditions $V_i(c(0))=e_i$. Then the normal frame is orthonormal along $c$, by invariance of $g$ and $V_i$ with respect to the flow of $X$.  Let $\varphi$ be an arbitrary smooth function on $[0,r]$ not vanishing everywhere and such that $\varphi(0)=\varphi(r)=0$. Define sections $W_i=\varphi V_i$ of $\VV$ along $c$. Then by Lemma \ref{lem-index1} we have $I(W_i,W_i)>0$, i.e.,
\[
\int_0^r((\varphi')^2g(V_i,V_i)-\varphi^2g(KV_i,V_i))dt=\int_0^r((\varphi')^2-\varphi^2g(KV_i,V_i))dt > 0,
\]
where prime denotes derivative with respect to $t$. Summing up with respect to $i$ we obtain
\[
\int_0^r(m(\varphi')^2-\tr(K)\varphi^2)dt>0.
\]
Since $\varphi$ was any nontrivial smooth function vanishing at $0$ and $r$ we can take $\varphi=\sin(at)$, where $a=\pi/r$. The above inequality and  $\tr K\ge \kappa$ imply 
\[\int_0^rma^2\cos^2(at)dt-\int_0^r\kappa\sin^2(at))dt>0.
\]
Taking into account that the integrals of $\cos^2(aT)$ and $\sin^2(at)$ over the period $[0,r]$ are equal we deduce that $ma^2>\kappa$. This implies that $r<\pi\sqrt{m/\kappa}=T^*$ and contradicts our supposition that $r>T^*$. This contradiction implies the assertion of the theorem.
\QED\medskip

The following lemma, proved below, can also be deduced from results in Chapter 11 of \cite{Hartman}.

\begin{lemma}\label{lem-index2}
Let $(P,Q)$ be $m\times m$ matrices defined on the interval $[0,T)$ (possibly infinite) which satisfy the differential equations
\[
\dot P=Q,\ \ \dot Q=-KP,
\]
and the initial conditions $P(0)=0$, $Q(0)=Id$, where $K$ is a symmetric matrix depending continuously on $t$. If $P$ is invertible on the interval $(0,r]$, for some $r<T$, then for any $C^1$-smooth $w:[0,r]\to \R^m$ satisfying $w(0)=0=w(r)$ we have
\begin{equation}\label{I3}
\int_0^r(\la \dot w,\dot w\ra-\la Kw,w\ra)dt\ge 0. 
\end{equation}
Equality takes place iff $w\equiv 0$ (the angle brackets denote standard scalar product in $\R^m$).
\end{lemma}

\paragraph{\bf Proof.}
Since $P$ is invertible on $I=(0,r]$ we may define $u$ by $w=Pu$ and then $u(0)=0=u(r)$ and $\dot w=\dot Pu+P\dot u$. We can decompose $L=\la \dot w,\dot w\ra-\la Kw,w\ra$ as
\[
L=\la P \dot u ,P\dot u\ra  +L_1+L_2
\]
where
\[
L_1=\la \dot P u,\dot P u\ra + \la \dot P u,P \dot u\ra, \quad L_2=\la P \dot u,\dot P u\ra -\la KPu,Pu\ra. 
\] 
We will prove that $\int_I(L_1+L_2)dt=0$. This equality implies the inequality 
\[
\int_IL\,dt=\int_0^r\la P\dot u,P\dot u\ra dt\ge 0
\]
and, by invertibility of $P$ and $u(0)=0$, we will have $\int_0^rLdt=0$ iff $u\equiv 0$ iff $w\equiv 0$. This will conclude the proof.

To prove $\int_I(L_1+L_2)dt=0$ note that $\dot P=Q$ implies
\begin{equation}\label{L1}
\int_I L_1dt=\int_I\la Qu,\dot{(Pu)}\ra dt=-\int_I \la \dot Qu,Pu\ra+\la Q\dot u,Pu\ra dt,
\end{equation}
where the latter equality follows from  integration by parts and $u(0)=0=u(r)$.
Using also $\dot Q=-KP$ we get
\begin{equation}\label{L2}
\int_I L_2dt=\int_I \la Qu,P\dot u\ra+\la \dot Qu,Pu\ra dt.
\end{equation}
We claim that $\la Qu,P\dot u\ra=\la Q\dot u,Pu\ra$ and consequently
\begin{equation}\label{L3}
\int_I \la Qu,P\dot u\ra-\la Q\dot u,Pu\ra dt=0. 
\end{equation}
This follows from 
\[
\la Qu,P\dot u\ra=\la P^TQu,\dot u\ra,\ \ \la Q\dot u,Pu\ra=\la P^TQ\dot u,u\ra
\]
and the symmetry of $P^TQ$ which is a consequence of
\[
\frac{d}{dt}P^TQ=\dot P^TQ+P^T\dot Q=Q^TQ-P^TKP.
\]
Here the right hand side is symmetric, by the assumption $K=K^T$, and the initial condition $P(0)=0$ implies symmetric $(P^TQ)(0)=0$, thus the solution $(P^TQ)(t)$ is symmetric.  From \eqref{L1}, \eqref{L2} and \eqref{L3} we see that $\int_I(L_1+L_2)dt=0$, which concludes the proof.
\QED
\medskip

\paragraph{\bf Proof of Lemma \ref{lem-index1}.}
Proposition 3.10 allows us to reduce the lemma to Lemma \ref{lem-index2} by the following argument. Let $\hat J=(J_1,\dots,J_m)$ be a tuple of Jacobi fields satisfying initial conditions $J_i(0)=0$ and $D_XJ_i(0)=e_i$, where $e_1,\dots,e_m$ is a basis in $\VV(c(0))$. Let $V_1,\dots,V_m$ be a normal frame along $c$ orthonormal with respect to $g$ which exists by Proposition 2.16. Expressing the Jacobi fields and the Jacobi frame in the normal frame used as a basis at $c(t)$ for all $t$ we can replace $\hat J(c(t))$ by its matrix, denoted $P(t)$. Then $D_X\hat J(c(t))$ corresponds to $\dot P$ and $D_X^2\hat J(c(t))$ corresponds to $\ddot P(t)$. Moreover, the Jacobi equation $D_X^2\hat J+K\hat J=0$ is equivalent to the equation $\ddot P+KP=0$ when the curvature operator is replaced with its matrix in the normal frame denoted again by $K$. This equation can be replaced by the system 
\[
\dot P=Q,\ \ \dot Q=-KP.
\]
By Proposition 3.10 a point $t>0$ is a conjugate point for the  dynamic pair $(X,\VV)$ iff there is a nontrivial Jacobi field $J$ such that $J(c(0))=0$ and $J(c(t))=0$. This is equivalent to the fact that the above Jacobi tuple $\hat J$ is linearly dependent at $c(t)$ and, in consequence, that the matrix $P$ is not invertible at $t$. In this way we see that the assumptions of the lemma imply the assumptions of Lemma \ref{lem-index2}. Namely, nonexistence of conjugate points in $[0,r]$ for the dynamic pair implies that the matrix $P$ defined here is invertible on $[0,r]$. Using Lemma \ref{lem-index2} we conclude that the integral functional in the lemma has property \eqref{I3}. 

To conclude we should show that the property \eqref{I3} in Lemma \ref{lem-index2} implies the similar statement \eqref{I2} on the functional \eqref{I1} in Lemma \ref{lem-index2}. Writing $V(c(t))=\sum_iw_i(t)V_i(c(t))$ and defining $w=(w_1,\dots,w_m)$ we have $w(0)=0=w(r)$, as $V(0)=0=V(r)$. Since the normal frame was chosen orthonormal with respect to $g$ we have
\[
I(V,V)=\int_0^r(g(D_XV,D_XV)-g(KV,V))dt= \int_0^r(\la \dot w,\dot w\ra-\la Kw,w\ra)dt  
\]
and the inequality \eqref{I3} in Lemma \ref{lem-index2} implies the inequality \eqref{I2} in Lemma \ref{lem-index1}. This completes the proof.
\QED\medskip

\paragraph{\bf Proof of Theorem \ref{thm3}.}
We first prove statement (ii). Note that existence along $c$ of a one dimensional subspace of $\VV$ parallel with respect to $D_X$ means that there exists a nonvanishing section $V\in\VV$ along $c$ which is normal, i.e., satisfies $D_XV=0$. Existence of $r$ such subspaces means that there are $r$ normal, linearly independent along $c$ sections  of $\VV$. These sections can be completed to a normal frame in $\VV$ along $c$ since such a frame is uniquely defined by a choice of a frame at the initial point $c(0)$ of $c$, by Proposition 2.6. The normal curvature matrix $K_t$ along $c$ has $r$ one dimensional diagonal blocks with the eigenvalues $\lambda_i(t)$. The Jacobi equation in Definition 3.8 written in this frame has $r$ single second order equations of the form \eqref{SJ} and an independent system of $m-r$ equations. Taking zero initial conditions for the values of all $m$ variables and for all their derivatives at $t=0$ with the exception of $\dot y_i(0)\not=0$ for a single $1\le i \le k$ we obtain $r$ linearly independent (over $R$) solutions of the Jacobi equation. If $k$ of these solutions vanish at some $t^*\in (0,T]$ then Proposition 3.10 implies that $t^*$ is a conjugate time of multiplicity at least $k$. This proves statement (ii).

Statement (i) follows from statement (ii) and the Sturm comparison theorem. Namely, the theorem implies that if $\lambda(t)$ and the constant $\kappa$ satisfy \eqref{LA} then a nontrivial solution on $[0,T]$ of the scalar equation $\ddot y=\lambda(t)y$ satisfying $y(0)=0$ has a zero between two consecutive zeros on $[0,T]$ of any nontrivial solution of the equation and $\ddot y'=\kappa y'$. Taking the solution $y'(t)=sin(\sqrt{\kappa}\,t)$ we see that $y(t)$ has a zero in $(0,\pi/\sqrt{\kappa}\,]$. 
\QED

\section{Semi-Hamiltonian systems}\label{s-semi-Hamiltonian}

Classical results on conjugate points assume presence of a Riemannian or Finslerian metric on a manifold $N$ and the vector field $X$ is the geodesic spray on the tangent bundle $M=TN$ while the vertical distribution $\VV$ is the distribution spanned by the vector fields tangent to the fibers of $TN\to N$. 
Another case where a canonical metric is given a priori is the case of $X$ being a Hamiltonian vector field on the cotangent bundle $T^*N$ and the distribution $\VV$ is tangent to the fibers of $T^*N\to N$. In a more general setting such a metric is given if and $\VV$ an integrable Lagrangian distribution on a symplectic manifold $(M,\sigma)$, cf. \cite{Agr,Chittaro}. The metric on $\VV$ is then given by $g(Y,Z)=\sigma([X,Y],Z)$. In our results below the closedness property of $\sigma$ and integrability of $\VV$ is not needed.

\begin{definition}\label{def-semi-Hamiltonian1}
Let $\sigma$ be a non-degenerate, skew-symmetric 2-form on a manifold $M$ of even dimension. 
We will call a vector field $X$ on $M$ \emph{semi-Hamiltonian} with respect to $\sigma$ if its flow preserves $\sigma$, i.e. $L_X\sigma=0$, where $L_X\sigma$ is the Lie derivative of $\sigma$ with respect to $X$. 

More generally, if $\DD$ is a distribution on $M$ of even rank and $\sigma$ is a non-degenerate, skew-symmetric 2-form on $\DD$ then $X$ will be called \emph{semi-Hamiltonian} (with respect to $(\DD,\sigma)$) if $\DD$ and $\sigma$ are invariant under the flow of $X$: 
\[
[X,\DD]\subset \DD, \eqno{(I)}
\]
\[
X(\sigma(Y,Z))=\sigma([X,Y],Z)+\sigma(Y,[X,Z]), \quad \forall\ Y,Z\in \DD. \eqno{(S)}
\]
\end{definition}
Clearly, if $\sigma$ is closed and $\DD=TM$ then we get the classical definition of Hamiltonian vector field.
If $\rk\DD=2m$ then we will call a rank $m$ sub-distribution $\VV\subset \DD$ \emph{Lagrangian} (with respect to $\sigma$) if the field of subspaces $\VV(x)\subset\DD(x)$ satisfies
\[
\sigma_x(v,w)=0, \quad \forall\ \  v,w\in\VV(x), \quad x\in M. \eqno{(L)}
\]

\begin{definition}\label{def-semi-Hamiltonian2}
A quadruple $(X,\VV,\DD,\sigma)$ where $X$ is semi-Hamiltonian with respect to $(\DD,\sigma)$, the sub-distribution $\VV\subset\DD$ is Lagrangian with respect to $\sigma$, and the regularity conditions (R1), (R2) hold, will be  called \emph{regular semi-Hamiltonian quadruple}.
\end{definition}

Above and furthermore all data are assumed to be $C^\infty$-smooth. The following systems, when satisfy (R1) and (R2), provide examples of regular semi-Hamiltonian quadruples.
\medskip

\begin{example} {\bf Time dependent Hamiltonian vector field.} 
Consider a time-varying vector field $\tilde X(t,x)$ on a symplectic manifold $(N,\tilde\sigma)$ and take:
\[
M=\R\times N, \quad  X=\frac{\partial}{\partial t}+ \tilde X, \quad \DD=\ker dt,
\]
where $t$ is the canonical coordinate on $\R$, extended to $M$ with the canonical projection $\R\times N\to \R$, the vector field $\tilde X$ is tangent
to the fibers $\{t\}\times N$ and identical with the original $\tilde X$ on $N$. The 2-form $\sigma$ is uniquely defined on $M$ by the form $\tilde \sigma$ and
the requirements: $\sigma\vert_{\{t\}\times N}=\tilde\sigma$ and $\sigma(\partial/\partial t,\cdot)=0$. Take $\VV(x)\subset\DD(x)$ to be a Lagrangian
subspace of $\tilde\sigma$ in $\DD(x)$, for any $x$, so that $x\mapsto \VV(x)$ is a smooth distribution. If $\VV+[X,\VV]=\DD$ then so defined quadruple $(X,\VV,\DD,\sigma)$ is regular semi-Hamiltonian.
\end{example}

\begin{example} {\bf Contact vector field on a contact manifold.} 
Let $(M,\omega)$ be a manifold of dimension $2m+1$ with a given 1-form $\omega$ on $M$ satisfying $\omega\wedge(d\omega)^m(x)\ne0$ for all $x$, called contact form. The contact form defines the rank $2m$ distribution $\DD=\ker\omega$ called contact distribution. 
Take $X$ to be a nonvanishing contact vector field on $M$, i.e., a nonvanishing vector field satisfying $[X,\DD]\subset\DD$, and let 
\[
\sigma=d\omega\vert_\DD.
\]
Finally, consider a Lagrangian sub-distribution $\VV\subset\DD$, $\rk\VV=m$, satisfying by definition $\sigma(Y,Z)=0$ for $Y,Z\in \VV$. Assume that the triple $(X,\VV,\DD)$ satisfies the regularity condition (R2). Then the curvature matrix $K$ and covariant differentiation operator $D_X$ on sections of $\DD$ are well defined (Sections \ref{ss-Curvature} and \ref{ss-Splitting}). Together with $\sigma$ we obtain the regular semi-Hamiltonian quadruple $(X,\VV,\DD,\sigma)$ on $M$. A special case is when $X$ is the Reeb vector field canonically defined by $\omega$ using formulas $\omega(X)=1$, $d\omega(X,\cdot)=0$. 
\end{example}

Any semi-Hamiltonian quadruple defines a bilinear form on $\VV$ by the formula
\begin{equation}\label{eq-metric}
g(V,W)=\sigma([X,V],W), \quad V,W\in\VV,
\end{equation}
which is tensorial since $\VV$ is Lagrangian. A regular semi-Hamiltonian quadruple $(X,\VV,\DD,\sigma)$ has the following properties.

\begin{proposition}\label{p8.1}
If $V_1,\dots,V_m$ form a normal frame in $\VV$ defined on an open subset $U\subset M$ then the following statements hold on $U$.

\noindent (i) The sub-distribution $\HH=\spn\{[X,V_1],\dots,[X,V_m]\}$ of $\DD$ is Lagrangian with respect to $\sigma$.

\noindent (ii) The matrix $(g_{ij})$ given by
\[
g_{ij}(x)=\sigma([X,V_i],V_j)(x), \quad i,j=1,\dots,m,
\]
is symmetric, nondegenerate, and defines a pseudo-Riemannian metric $g$ on $\VV$.

\noindent (iii) The metric $g$ is invariant under covariant differentiation along $X$, i.e.,
\begin{equation}
X(g(Y,Z))=g(D_XY,Z)+g(Y,D_XZ).
\end{equation}

\noindent (iv) The matrix $K$ (Definition \ref{d2.1}) is symmetric and defines a selfadjoint (relative to $g$) operator $K\colon\VV\to\VV$. 
\end{proposition}

\begin{remark}
It can be deduced (e.g. exploiting structural equations from \cite{Agr}) that the curvature operator $K$ used here coincides, in a Hamiltonian setting, with a curvature introduced in the work of Agrachev and coworkers, cf. \cite{AgrGam1,AgrGam2,Agr, AgrZel}, that starts from an optimality problem and goes through lifting of the system data along an extremal to a curve in a Lagrangian grassmannian. Via symplectic geometry arguments it is shown that the curve contains essential information on the optimality problem and the curvature is one of its invariants. So defined curvature is more geometric but difficult to compute, in general.  The curvature (curvatures) defined in \cite{JK} and here was (were) introduced with the motivation of understanding better the dynamic properties of trajectories of systems of ODEs and of classes of control systems, without stating any optimality problem. As stated earlier, in the case of systems of ODEs of second order such curvature is classic.
\end{remark}

Given the metric $g$ and the curvature operator $K$, we can define the directional curvature $k(v)=g(Kv,v)/g(v,v)$, for $v\in \VV$ such that $g(v,v)\ne0$, as in Section \ref{ss-Results}. Theorems \ref{thm1} and \ref{thm2} have the following counterparts in the semi-Hamiltonian setting.

\begin{theorem}\label{thmH1}
Consider a regular semi-Hamiltonian quadruple $(X,\VV,\DD,\sigma)$ on $M$ and let $c:[0,T]\to M$ be an integral curve of $X$. We assume that:

\smallskip
\noindent (A) the metric \eqref{eq-metric} is strictly definite (positive or negative) at all points of $c$;

\noindent (B) the directional curvature satisfies $k(v)\le \lambda$ for all $v\in \VV(x)\setminus\{0\}$ and all $x$ in $c$, where $\lambda\in\R$ is a constant.

Then there are no conjugate times in $(0,T]$, if $\lambda\le 0$, and there are no conjugate times in the interval $(0,t_c)$,
$t_c=\min\{T,\frac{\pi}{\sqrt{\lambda}}\}$, if $\lambda>0$.
\end{theorem}
\paragraph{\bf Proof.}
The proof is the same as the proof of Theorem \ref{thm1}, with the following modification. Instead of using an "artificial" metric introduced there we use
the canonical metric $g(Y,Z)=\sigma([X,Y],Z)$ on $\VV$ introduced above. Assumption (A) guarantees that it is definite and we may assume that it is
positive definite (otherwise we replace it with $-g$). The metric is invariant under covariant differentiation with respect to $D_X$, as Proposition
\ref{p8.1} says. Introducing the function $\varphi=\sqrt{g(V,V)}$ along a geodesic, where $V$ is a vertical Jacobi field, we can proceed in exactly the
same way as in the proof of Theorem \ref{thm1}, thus proving Theorem \ref{thmH1}.
\QED

\begin{theorem}\label{thmH2}
Consider a regular semi-Hamiltonian quadruple $(X,\VV,\DD,\sigma)$ on $M$ and let $c:[0,T)\to M$ be an integral curve of $X$. We assume that (A) is
satisfied and $ \tr\,K \geq \kappa>0$, along $c$.
Then there is a conjugate point in $(0,T^*]$ if  $T^*=\pi\frac{\rk\VV}{\sqrt{\kappa}}<T$.
\end{theorem}
\paragraph{\bf Proof.}
The corollary follows from Theorem \ref{thm2} and Proposition \ref{p8.1}. Namely, for the metric $g$ in Theorem \ref{thm2} we take the metric from Proposition \ref{p8.1},
if it is positive definite (otherwise we take $-g$). By the proposition the curvature $K$ is symmetric with respect to such $g$ and the assertion follows from Theorem \ref{thm2}.
\QED\medskip

\paragraph{\bf Proof of Proposition \ref{p8.1}.}
To prove statement (ii) take $V_i,V_j\in\VV$. Since $\VV$ is Lagrangian, we have $\sigma(V_i,V_j)=0$. Lie differentiating this equality with respect to $X$ and using condition (S) gives
$$
0=\sigma([X,V_i],V_j)+\sigma(V_i,[X,V_j])
$$
and proves the symmetry of $g$ in (ii), since $\sigma$ is antisymmetric. 

Differentiating  twice gives
$$
0=\sigma(\ad^2_XV_i,V_j)+2\sigma([X,V_i],[X,V_j])+\sigma(V_i,\ad^2_XV_j).
$$
The side terms are zero which follows from \eqref{eq13} and the fact that $\VV$ is Lagrangian. Thus
$\sigma([X,V_i],[X,V_j])=0$, which shows (i).

Differentiating three times yields
$$
0 = \sigma([X,\ad_X^2V_i],V_j) + 3\sigma(\ad_X^2V_i,[X,V_j]) + 3\sigma([X,V_i],\ad_X^2V_j) + \sigma(V_i,[X,\ad_X^2V_j])
$$
and, applying \eqref{eq13} and the summation convention,
\begin{eqnarray*}
0 & = & \sigma([X,K_i^sV_s],V_j) + 3\sigma(K_i^sV_s,[X,V_j])\\
   & &  +\ 3\sigma([X,V_i],K_j^sV_s) + \sigma(V_i,[X,K^s_jV_s])\\
  & = &  \sigma(K_i^s[X,V_s],V_j) + 3\sigma(K_i^sV_s,[X,V_j])\\
   & & +\ 3\sigma([X,V_i],K_j^sV_s) + \sigma(V_i,K^s_j[X,V_s]),
\end{eqnarray*}
where in the second equality we use  the fact that $\sigma(X(K^s_i)V_s,V_j)=0$, as $\VV$ is Lagrangian.
Using antisymmetry of $\sigma$ and the definition of $g_{ij}$ we get
$$
0=K_i^sg_{sj}-3K_i^sg_{js}+3K_j^sg_{is}-K_j^sg_{si}=-2K_i^sg_{sj}+2K_j^sg_{is},
$$
since $g$ is symmetric. This proves (iv).

Since $\VV=\spn\{V_1,\dots,V_m\}$ and $\HH^1=\spn\{[X,V_1],\dots,[X,V_m]\}$ are transversal Lagrangian subspaces
in $TN$ and $\sigma$ is non-degenerate on $TN$, it follows that the matrix $g$ in (ii) is non-degenerate and defines a pseudo-Riemannian metric on $\VV$.

Finally, for proving statement (iii) it is enough to show that the coefficients $g_{ij}$ are constant along $c$. Using invariance of $\sigma$ with respect to $X$ (condition (S)) we find that $X(g_{ij})=\sigma([X[X,V_i]],V_j)+\sigma([X,V_i],[X,V_j])=0$ since $V_i$ and $V_j$ are normal, thus $[X[X,V_i]]$ and $V_j$ are in the Lagrange distribution $\VV$ and $[X,V_i]$, $[X,V_j]$ are in the Lagrange distribution $\HH$. 
\QED

\begin{remark}
The classical example of a Hamiltonian system on the phase space $M=T^*Q$ fits into above framework. Namely, given a function $H:T^*Q\to \R$ (Hamiltonian) on the cotangent bundle of a (configuration) manifold $Q$ and the canonical symplectic structure $\sigma$ on $T^*Q$, take $X=\vec H$ the Hamiltonian vector field on $T^*Q$ defined by $H$  and let $\VV$ be the vertical distribution on the bundle $\pi:T^*Q\to Q$. The tuple $(\vec H,\VV,\DD,\sigma)$, where $\DD=T(T^*Q)$, forms a regular Hamiltonian quadruple in the region where the Hamiltonian $H$ is regular, i.e., its derivative along the fibers is is nonzero and the square matrix of its second derivatives along the fibers is nondegenerate. The curvature operator $K$ is then the classical Jacobi endomorphism of the Hamiltonian equations defined by $\vec H$. 

A regular semi-Hamiltonian quadruple is obtained if we restrict the considerations to a level submanifold $M_c=\{H=c\}\subset T^*Q$ of the Hamiltonian $H$ and take the vector field $X$ equal to $\vec H$ restricted to $M_c$, $X=\vec H\vert_{M_c}$. The distribution $\VV$ on $M_c$ is defined as the vertical distribution of the cotangent bundle intersected with the tangent space to $M_c$: $\VV(x)=T(T^*_{\pi(x)}Q)\cap T_xM_c$, for $x\in M_c$. Then, if $\dim Q=m+1$, we have $\dim \VV(x)=m$, $\dim M_c=2m+1$ and assuming regularity of such pair $(X,\VV)$ makes sense (typical examples are regular). In this case all statements of Proposition \ref{p8.1} hold true for the canonical symplectic form $\sigma$ on $T^*Q$ replaced with the skew-symmetric form $\hat\sigma=\sigma\vert_{M_c}$ restricted to the corresponding distribution $\DD$ on $M_c$. 
\end{remark}

\section{Special cases}\label{s-Examples}

\begin{example}\label{ex-spray}{\bf Geodesic spray.}
Let $(N,g)$ be a Riemannian or pseudo-Riemannian manifold. Take $M=TN$ to be the tangent bundle of $N$, $\VV$ the distribution tangent to the fibers of
$\pi:TN\to N$, and $X$ the geodesic spray corresponding to metric $g$. In local coordinates $X=\sum_iy^i\partial_{x^i}+\sum_{i,j,k}\Gamma_{jk}^i(x)
y^jy^k\partial_{y^i}$ and $\VV=\spn\{\partial_{y^i}\}$. Clearly, $(X,\VV)$ is a regular pair. It can be shown (see \cite[Section 2.3]{JK}) that $K(x,y)$, $(x,y)\in M=TN$ satisfies the equality $K(x,y)v=R_x(v,y)y$, where $R_x$ is the Riemann curvature tensor at $x\in N$ and $v$, $y$ are considered as vectors in $T_xN$. The same holds in the case of Finsler manifolds and, more generally, spray manifolds in which case the coefficients $\Gamma_{jk}^i$ are functions of $(x,y)$, instead of $x$ only. Our theorems imply then the corresponding versions of Cartan--Hadamard and Bonet--Myers
theorems on Riemannian, Finsler and spray manifolds, respectively.
\end{example}

\begin{example} {\bf Control-affine systems.}
In Geometric Control Theory the most often considered class of non-linear systems are control-affine systems,
\begin{equation}\label{ca-system}
\dot x=f(x)+\sum_{i=0}^mu_ig_i(x),
\end{equation}
where $f,g_1,\dots,g_m$ are vector fields on a manifold $M$, $x\in M$ is the state of the system, and  $u_1,\dots,u_m$ are components of the control.
Such a system defines a dynamic pair,
\[
X=f,\quad \VV=\spn\{g_1,\dots,g_m\}.
\]
The regularity conditions (R1), (R2) and the invariance condition (I) are often satisfied for systems met in applications. Assuming that $\dim M=2m$,
conditions (R1), (R2) mean that $f(x)\not=0$, $\spn\{g_1,\dots,g_m,[f,g_1],\dots,[f,g_m]\}=TM=\DD$.

These assumptions are satisfied for so called fully actuated mechanical systems, met in Robotics. For such a system $M=TQ$ is the velocity phase space,
where $Q$ is the configuration manifold. The system is usually described (see \cite{BL} and also \cite{NR,RR}) by the Euler-Lagrange equations with external forces, corresponding to a Lagrangian
\[
L(q,\dot q)=\frac{1}{2}g(\dot q,\dot q)-V(q),
\]
where $q\in Q$ is the generalized position, the first term in $L$ is the kinetic energy (mathematically, a Riemann metric $g$ on $Q$) and $V$ is the
potential energy. The controlled equations of motion can be written, in coordinates $q^i$, $v^i=\dot q^i$, as
\begin{eqnarray}\label{mc-system}
\dot q^i & = & v^i, \nonumber\\
\dot v^i & = & -\sum_{jk} \Gamma^i_{j,k}(q)v^jv^k + \sum_sg^{is}\left(F^P_s+F_{0,s}(q,v)+\sum_ju^jF_{j,s}(q,v)\right).
\end{eqnarray}
Here $\Gamma^i_{jk}$ are coefficients of the Levi-Civita connection of the metric $g$ and $(g^{ts})$ is the inverse matrix to the matrix of metric
coefficients $(g_{ij})$ of $g$. The components $F^P_s$, $F_{0,s}$ are coefficients of potential and external uncontrollable forces, respectively, while
$F_{j,s}$ are coefficients of a controlled external force $F=\sum_ju_jF_j$. The system is fully actuated if $F_1,\dots,F_m$ are linearly independent and
$m=\dim Q$. The terms on the right-hand-side not containing controls define a vector field $X$ on $M=TQ$, while the vector fields $\hat
F_j=\sum_iF^i_j\partial_{v^i}$ on $M$ define control vector fields $g_j$ in \eqref{ca-system} and span the distribution $\VV$. A system \eqref{mc-system}
is called \emph{fully actuated} if $\VV=\spn\{\partial_{v^1},\dots,\partial_{v^m}\}$. This is the same distribution as the one defined by a second order
ODE, thus we have the following

\medskip
{\bf Claim.} {\it The fully actuated system \eqref{mc-system} defines the same dynamic pair (with $\VV$ as above) as the one assigned to the second order ODE obtained from \eqref{mc-system} by choosing controls equal to zero, i.e., the ODE given by the vector field on $M=TQ$
$$
X=\sum_iv^i\partial_{q^i}+\left(-\sum_{jk} \Gamma^i_{j,k}(q)v^jv^k + \sum_sg^{is}\left(F^P_s+F_{0,s}(q,v)\right)\right)\partial_{v^i}.
$$}

Consider a special case where the metric coefficients are constant (the kinetic energy depends only on the components of the velocities and not on the
positions) and the uncontrolled external force is independent of the velocity. Then $\Gamma^i_{j,k}=0$ and
\[
X=\sum_iv^i\partial_{q^i}+ \sum_sg^{is}(-P_{q^s}+F_{0,s}(q))\partial_{v^i},
\]
where $P_{q^s}$ are partial derivatives of the potential energy $P$. Choosing $V_j=\partial_{v^j}$ we see that $[X,V_j]=-\partial_{q^j}$ and
$[X,[X,V_j]]=-\sum K^i_j V_j$, where
\[
K^i_j=\sum_s g^{is}\left(-\frac{\partial^2 P}{\partial q^s\partial q^j}+\frac{\partial F_{0,j}}{\partial q^s}\right).
\]
This means that the frame $V_j=\partial_{v^j}$ is a normal frame of $\VV$ and $K=(K^i_j)$ is a normal curvature matrix, which can be used in Theorems
\ref{thm1} and \ref{thm2}.

Note that the external force $F_0$ can be split into two components $F_0=\bar F_0+\hat F_0$ so that the matrix of partial derivatives of the components
$(\partial \bar F_{0,s}/\partial q^t)$ is symmetric and $(\partial \hat F_{0,s}/\partial q^t)$ is antisymmetric. The first part can be included into the
potential force thus, without loosing generality, we can assume in the above formula for $K^i_j$ that the matrix $\frac{\partial F_{0,j}}{\partial q^s}$
is antisymmetric and represents a non-potential force.

Even if the antisymmetric part is zero, the above curvature matrix is nonsymmetric. This is due to the fact that the canonical Euclidean metric is not
natural in this case. The natural metric is given by the metric coefficients $(g_{ij})$, which are assumed constant. The corresponding orthonormal frame
$\hat V'=(V'_1,\dots,\tilde V'_m)$ is also a normal frame of $\VV$, related to $\hat V=(V_1,\dots,V_m)$ via a transformation $\hat V'=\hat VG$, with a
symmetric matrix $G$ such that $GG=(g^{ij})$. In this frame the curvature matrix, due to Proposition \ref{p2.1}, is $K'=G^{-1}\hat KG$ where
\[
\hat K=(\hat K^i_j), \qquad \hat K^i_j=-\frac{\partial^2 P}{\partial q^s\partial q^j}+\frac{\partial F_{0,j}}{\partial q^s}.
\]
\end{example}

\begin{example}\label{ex-SODE}{\bf A class of ordinary differential equations.}
Consider a pair of second-order ordinary differential equations:
\[
\ddot x_i=F_i(t,x,\dot x), \qquad i=1,2.
\]
If the corresponding normal curvature matrix $t\mapsto K(t)$ is diagonalizable, then, depending on the signs of the eigenvalues, Theorems \ref{thm1} and \ref{thm2} can be utilized to obtain estimates for conjugate points. However, identifying a general class of such systems is a challenging task. Nonetheless, such systems have been encountered in a different context in \cite{KM}, where they were identified as those arising from the so-called "dancing construction" for path geometries (also discussed in \cite{BLN,D}). The following specific example illustrating this concept was found in \cite[Example 4.17]{KM}
\[
\ddot x_1= F(t, x_1, \dot x_1),\qquad \ddot x_2=\frac{\dot x_2}{\dot x_1-x_2}(F(t,x_1,\dot x_1)-2\dot x_2)
\]
where $F$ is arbitrary smooth function of 3 variables. Recall that $\VV=\spn\{\partial_{y_1},\partial_{y_2}\}$, where $y_i$ corresponds to $\dot x_i$. The curvature operator is given in the natural basis $\partial_{y_1},\partial_{y_2}$ by
\[
K=
\left(
\begin{array}{cc}
\chi_1&0\\
-\frac{y_2(\chi_1+\chi_2)}{y_1-x_2}&\chi_2\\
\end{array}
\right)
\]
where $\chi_1=-\partial_{x_1}F+\frac{1}{2}X(\partial_{y_1}F)-\frac{1}{4}(\partial_{y_1}F)^2$ is the curvature for the single ODE $\ddot x_1=F(t,x_1,\dot x_1)$, and $\chi_2=\frac{1}{2}\frac{X(F)}{y_1-x_2}+\frac{3}{4}\frac{F(2y_2-F)}{(y_1-x_2)^2}$ (recall that $X$ denotes here the total derivative, as usual in the context of ODEs). Clearly, $K$ is not-symmetric but it is diagonalizable and the corresponding eigen-vectors are $V_1=\partial_{y_2}$ and $V_2=\partial_{y_1}+\frac{y_2}{y_1-x_2}\partial_{y_2}$. The two directions are preserved by $D_X$. Indeed, applying Proposition \ref{p-horizontal} and Definition \ref{d4.3} with $H_1$ given in Section \ref{ss-ODEs} we compute that
\[
D_X V_1=-\frac{1}{2}\left(\frac{F-4y_2}{y_1-x_2}\right)V_1\qquad\mathrm{and}\qquad D_XV_2=-\frac{1}{2}(\partial_{y_1}F)V_2.
\]
Hence, it follows that $V_i$, $i=1,2$, can be rescaled to get a normal frame, with $\operatorname{diag}(\chi_1,\chi_2)$ as the normal curvature matrix.
Note that the two eigen-values $\chi_1$ and $\chi_2$ can be of different signs, depending on the particular choice of function $F$. In particular $K=0$ for $F=0$, which corresponds to a model of simultaneous moves (a dance) of points in the projective surfaces $\R P^2$ and the dual projective surface $(\R P^2)^*$, c.f. \cite{BLN,D,KM}. 
\end{example}

\begin{example}\label{s-Example}{\bf Non-symmetric perturbations.}
Consider the family of systems
\begin{equation}\label{example1}
\ddot x=-x-\e y,\quad \ddot y=-y+\e x
\end{equation}
parametrized by $\e\in\R$.
For $\e=0$ the system has double conjugate times $t_k=k\pi$,  $k=1,2,\dots$, and the curvature matrix is the identity matrix $K=I$.
We will show that for arbitrary $N>0$ the first $N$ of them disappear after a small perturbation with $\e\not=0$, $|\e|\le \delta$, if $\delta$ is small enough (note that then $K$ is nonsymmetric).

In order to compute solutions for arbitrary $\e$ it is convenient to replace the vector $(x,y)$ by the complex number $z=x+iy$. Then \eqref{example1} can be written as
$\ddot z=-(1+i\e)z$
or as a system of first order equations
\begin{eqnarray*}
\dot z & = & w,\\ \dot w & = & -(1+i\e)z.
\end{eqnarray*}
Taking initial conditions $z(0,\e)=0$, $w(0,\e)=w_0$ we see that the conjugate times will be given by the values of $t>0$ such that $z(t,\e)=0$.
Introducing the row vector $C=(1,0)$, the column vector $W=(0,w_0)$ and the system matrix
\[
A=\begin{pmatrix} 0 & 1 \\ -1-i\e & 0 \end{pmatrix}
\]
we see that $z(t,\e)=C\exp(tA) W$. Since $A^2=-(1+i\e)I$, $A^{2n}=(-1)^n(1+i\e)^nI$, $A^{2n+1}=(-1)^n(1+i\e)^nA$, $CA^{2n}W=0$, $CA^{2n+1}W=(-1)^n(1+i\e)^nw_0$ we find that
\begin{eqnarray*}
z(t,\e) & = & \sum_{n\ge 0}\frac{t^n}{n!}CA^nW\ =\ \sum_{n\ge 0}\frac{t^{2n+1}}{(2n+1)!}CA^{2n+1}W  =  \sum_{n\ge 0}\frac{t^{2n+1}}{(2n+1)!}(-1)^n(1+i\e)^nw_0\\
        & = & \sum_{n\ge 0}\sum_{0\le k\le n}(-1)^n\frac{t^{2n+1}}{(2n+1)!}\binom{n}{k}(i\e)^kw_0  =  \sum_{k\ge 0}\ \sum_{n\ge k}(i\e)^k(-1)^n\binom{n}{k}\frac{t^{2n+1}}{(2n+1)!}w_0\\
        & = & \sum_{k\ge 0}(i\e)^kf_k(t)w_0,
\end{eqnarray*}
with suitable real analytic functions $f_k(t)$.
Note that the coefficients $a_{nk}$ of the middle double sum can be majorized by $t|w_0|(t^2)^n\e^k(n!k!)^{-1}$, thus the summation can be done in arbitrary order. We deduce that
\[
z(t,\e)= (F_{re}(t,\e)+i\e F_{im}(t,\e))w_0,
\]
where
\[
F_{re}(t,\e)=f_0(t)+\e^2\tilde F_{re}(t,\e),  \quad F_{im}(t,\e)=f_1(t)+\e^2\tilde F_{im}(t,\e)
\]
and $\tilde F_{re}(t,\e)$, $\tilde F_{im}(t,\e)$  are real analytic functions.
Moreover, we have
\[
f_0(t)=\sum_{n\ge 0}(-1)^n\frac{t^{2n+1}}{(2n+1)!}=\sin(t),
\]
which means that the system with $\e=0$ has double conjugate times at
$t_k=k\pi$, $k=1,2,\dots$. We claim that they disappear for small enough $\e\not=0$.
More precisely,

\medskip
{\bf Claim 1.} {\it For any integer $N>0$ there exists $\delta>0$ such that for $\e\not=0$, $|\e|<\delta$ the system
does not have conjugate times in $(0,N\pi]$.}

\medskip

For proving this it is enough to show that the functions
$F_{re}(t,\e)$ and $F_{im}(t,\e)$ can not be simultaneously zero for $t\in (0,N\pi]$ and small enough $\e\not=0$.
We will use the following

\medskip
{\bf Claim 2.} {\it The function
\[
f_1(t)=\sum_{n\ge 1}(-1)^n \frac{n\,t^{2n+1}}{(2n+1)!} = - \frac{1}{2}\int_0^ts \sin(s) ds =  \frac{1}{2}(-\sin t+t\cos t)
\]
oscillates between its local maxima and minima taking values $f_1(k\pi)=\frac{1}{2}(-1)^k k\pi$ and it is strictly monotonous in the intervals $[k\pi,(k+1)\pi]$.}
\medskip

To show Claim 1 fix $N>0$ and denote $I_N=[0,N\pi]$.
It is enough to find $\delta>0$ such that for any fixed $\e$ with $|\e|<\delta$ the functions $F_{re}(t,\e)$ and $F_{im}(t,\e)$ have
no common zeros in $I_N$ or, equivalently, $G_\e(t)=|F_{re}(t,\e)|+|F_{im}(t,\e)|>0$ for $t\in I_N$.
Note that for $\e=0$ the function $G_0(t)=|f_0(t)|+|f_1(t)|$ does not vanish in $I_N$ as $f_0(t)=sin t$ has zeros at $t=k\pi$,
while $f_1(t)$ is nonzero at such points, by Claim 2. Thus
\[
m:=\min_{t\in I_N}\{\,|f_0(t)|+|f_1(t)|\,\}>0
\]
since $I_N$ is compact. Moreover, for $t\in I_N$ we have
\begin{eqnarray*}
G_\e(t) & = & |f_0(t)+\e^2\tilde F_{re}(t,\e)|+|f_1(t)+\e^2\tilde F_{im}(t,\e)|\ge \\
        & = & |f_0(t)|-\e^2|\tilde F_{re}(t,\e)| + |f_1(t)|-\e^2|\tilde F_{im}(t,\e)|\ge  m-\e^2(a+b),
\end{eqnarray*}
where $a=\max_{t\in I_N}\{|\tilde F_{re}(t)|\}>0$ and $b=\max_{t\in I_N}\{|\tilde F_{im}(t)|\}>0$.
It follows that $G_\e(t)>0$, if $\e^2<m(a+b)^{-1}$, which shows Claim 1.

It remains to prove Claim 2. The first equality is just the definition of $f_1$ and the last one follows by integration by parts.
The middle equality is a consequence of the calculation
\begin{eqnarray*}
f_1(t) & = & \frac{1}{2}\sum_{n\ge 1}(-1)^n \frac{1}{(2n-1)!}\frac{t^{2n+1}}{2n+1} = \frac{1}{2}\sum_{n\ge 1}(-1)^n \frac{1}{(2n-1)!}\int_0^ts^{2n}ds\\
       & = & \frac{1}{2}\int_0^ts\sum_{n\ge 1}(-1)^n \frac{1}{(2n-1)!}s^{2n-1}ds = - \frac{1}{2}\int_0^ts \sin(s) ds.
\end{eqnarray*}
\QED
\end{example}


\begin{thebibliography}{10}
\bibitem{Agr} A. Agrachev, \textit{The curvature and hyperbolicity of Hamiltonian systems}, Proceed. Steklov Math. Inst., vol. 256 (2007), 26--46.
\bibitem{AB} A. Agrachev, I. Beschastny, \textit{Jacobi fields in optimal control: Morse and Maslov indices}, Nonlinear Analysis 214 (2022), 112608.
\bibitem{AgrGam1} A. Agrachev, R. Gamkrelidze, \textit{Feedback--invariant optimal control theory and differential geometry, I. Regular extremals}, J. Dynamical and Control Systems, vol. 3 (1997), 343--389. 
\bibitem{AgrGam2} A. Agrachev, R. Gamkrelidze, \textit{Feedback--invariant optimal control theory and differential geometry, II. Jacobi curves for singular extremals}, J. Dynamical and Control Systems, vol. 4(4) (1998), 583--604. 
\bibitem{AgrZel}  A. Agrachev, I. Zelenko, \textit{Geometry of Jacobi curves II}, Journal of Dynamical and Control Systems volume 8, pages167--215 (2002)
\bibitem{BR} D. Barilari, L. Rizzi, \textit{Comparison theorems for conjugate points in sub-Riemannian geometry}, ESAIM: COCV 22(2) (2016), 439--472.
\bibitem{BLN} G. Bor, L. Lamoneda, P. Nurowski, \textit{The dancing metric, G2-symmetry and projective rolling},Trans. Amer, Math. Soc., vol. 370(6) (2018), 4433--4481.
\bibitem{BL} F. Bullo, A. Lewis, \textit{Geometric Control of Mechanical Systems}, Texts in Applied Mathematics (TAM, volume 49), Springer, 2005.
\bibitem{CGM} J.F. Carinena, I. Gheorghiu, E. Martínez, \textit{Jacobi fields for second-order differential equations on Lie algebroids}, In: Dynamical Systems and Differential Equations, Proceedings of the 10th AIMS International Conference, eds. M. de Leónet al. (Madrid, Spain, 2015), pp. 213--222.
\bibitem{CM} J.F. Carinena, E. Martinez, \textit{Generalized Jacobi equation and inverse problem in classical mechanics}, In: V.V. Dondonov and V. Manko (Eds), Integral Systems, Solid State Physics and Theory of Phase Transitions. Nova Science, 1992.
\bibitem{Chittaro} F.C. Chittaro, \textit{An estimate for the entropy of Hamiltonian flows}, J. Dynamical and Control Systems vol. 13 (2007), 55–67.
\bibitem{CMS} M. Crampin, E. Martinez and W. Sarlet, \textit{Linear connections for systems of second-order ordinary differential equations}, Ann. Inst. Henri Poincare, 65 (1996), 223--249.
\bibitem{CF} N. Caroff, H. Frankowska, \textit{Conjugate points and shocks in nonlinear optimal control}, Trans. Amer. Math. Soc. 348 (1996), 3133--3153 .
\bibitem{Cartan} E. Cartan, \textit{Observations sur le m\'emoir pr\'ec\'edent},  Math. Zeitschrift 37 (1933), 619-622.
\bibitem{Chern} S.-S. Chern, \textit{Sur la g\'eom\'etrie d'un syst\'eme d'\'equations diff\'erentialles du second ordre}, Bull. Sci. Math. 63 (1939), 206–212.
\bibitem{DMJ} J.Diaz, J. B. McLaughlin, R. Joyce, \textit{Sturm comparison theorems for ordinary and partial differential equations}, Bull. Amer. Math. Soc. 75 (1969), 335–339.
\bibitem{D} M. Dunajski, \textit{Twistor Theory of Dancing Paths}, SIGMA 18 (2022), 027, 13 pages.
\bibitem{Grossman} D. A. Grossman, \textit{Torsion-free path geometries and integrable second order ODE systems}, Selecta Math. (N.S.), 6(4):399--442, (2000).
\bibitem{HM} S. Hajd\'u, T. Mestdag, \textit{Conjugate points for systems of second-order ordinary differential equations}, Int. J. Geom. Methods Mod. Phys., 17(1) (2020).
\bibitem{Hartman} P. Hartman, \textit{Ordinary differential equations}, John Wiley and Sons, New York 1964.
\bibitem{JK} B. Jakubczyk, W. Kryński, \textit{Vector fields with distributions and invariants of ODE's}, J. Geom. Mech. 5(1) (2013) 85-129.
\bibitem{JP} M. Jerie, G.E. Prince, \textit{Jacobi fields and linear connections for arbitrary second-order ODEs}, J. Geom. Phys.(2002) 43(4), 351--370.
\bibitem{Kosambi} D.D. Kosambi, \textit{Parallelism and path-space},  Math. Zeitschrift 37 (1933), 608-618.
\bibitem{KM} W. Kry\'nski, O. Makhmali, \textit{A characterization of chains and dancing paths in dimension three}, arXiv:2303.08807 (2023).
\bibitem{NR} M. Nowicki, W. Respondek, \textit{A Mechanical Feedback Classification of Linear Mechanical Control Systems}, Special Issue ``Advances in Robot Motion and Control'', Appl. Sci. 2021, 11(22), 10669.
\bibitem{RR} S. Ricardo, W. Respondek, \textit{When is a control system mechanical? }, J. Geom. Mech., Volume 2, Issue 3: (2010) 265--302.
\bibitem{Sabau} S.V. Sabau, \textit{Some remarks on Jacobi stability}, Nonlinear Analysis 63 (2005) e143 – e153.
\end{thebibliography}
\end{document}